# MATRIX REPRESENTATIONS AND INDEPENDENCIES IN DIRECTED ACYCLIC GRAPHS

By Giovanni M. Marchetti[1] and Nanny Wermuth[2]

*University of Florence and Chalmers/Göteborgs Universitet*

For a directed acyclic graph, there are two known criteria to decide whether any specific conditional independence statement is implied for all distributions factorized according to the given graph. Both criteria are based on special types of path in graphs. They are called separation criteria because independence holds whenever the conditioning set is a separating set in a graph theoretical sense. We introduce and discuss an alternative approach using binary matrix representations of graphs in which zeros indicate independence statements. A matrix condition is shown to give a new path criterion for separation and to be equivalent to each of the previous two path criteria.

**1. Introduction.** We consider stepwise processes for generating joint distributions of random variables $Y_i$ for $i = 1, \ldots, d$, starting with the marginal density of $Y_d$, proceeding with the conditional density of $Y_{d-1}$ given $Y_d$, up to $Y_1$ given $Y_2, \ldots, Y_d$. The conditional densities are of arbitrary form but have the independence structure defined by an associated directed acyclic graph in $d$ nodes, in which node $i$ represents variable $Y_i$. Furthermore, an arrow starting at node $j > i$ and pointing to $i$, the *offspring* of $j$, indicates a non-vanishing conditional dependence of $Y_i$ on $Y_j$. Node $j$ is then called a parent node of $i$, the set of parent nodes is denoted by $\mathrm{par}_i \subseteq \{i+1, \ldots, d\}$, and the graph together with the complete ordering of the nodes as $V = (1, \ldots, d)$ is the *parent graph* $G^V_{\mathrm{par}}$.

Received December 2007.

[1]Supported in part by MIUR, Rome, under the Project PRIN 2005132307.

[2]Supported in part by the Swedish Research Society via the Gothenburg Stochastic Center and by the Swedish Strategic Fund via the Gothenburg Math. Modelling Center.

The joint density $f_V(y)$ of the $d \times 1$ random variable $Y$ factorizes as

$$(1.1) \qquad f_V(y) = \prod_{i=1}^{d} f_{i|\,\mathrm{par}_i}(y_i|y_{\mathrm{par}_i}),$$

where $f_{i|\,\mathrm{par}_i}(y_i|y_{\mathrm{par}_i}) = f_i(y_i)$ whenever $\mathrm{par}_i$ is empty.

For each node $j > i$, not a parent of $i$, the factorization (1.1) implies that the $ij$-arrow is missing in the graph and that $Y_i$ is conditionally independent of $Y_j$, given $Y_{\mathrm{par}_i}$. The defining list of independencies for $G_{\mathrm{par}}^V$, written in terms of nodes, is

$$(1.2) \qquad i \perp\!\!\!\perp j|\mathrm{par}_i \quad \text{for } j > i \text{ not a parent of } i \text{ and } i = 1, \dots, d-1.$$

In general, for any disjoint subsets $\alpha$, $\beta$ and $C$ of $V$ we denote $Y_\alpha$ conditionally independent of $Y_\beta$ given $Y_C$ by $\alpha \perp\!\!\!\perp \beta|C$. The parent graph is a special type of *independence graph*, that is, a graph for which each missing edge corresponds to an independence statement. The *independence structure captured by a graph* consists of the list of independencies defining the graph and of all other independences that derive from this list.

For instance, the graph $G_{\mathrm{par}}^V$ typically captures more independence statements than those given directly by the defining list (1.2). We take as an example the following graph representing a Markov chain in four nodes:

$$1 \prec\!\!-2 \prec\!\!-3 \prec\!\!-4.$$

The defining independence for node 1 is $1 \perp\!\!\!\perp \{3,4\}|2$, but for instance the additional independence statement $1 \perp\!\!\!\perp 4|3$ also holds. We say that a *distribution is generated over* $G_{\mathrm{par}}^V$ if its density $f_V(y)$ is obtained by the stepwise process described above so that it factorizes as in (1.1) and its set of independence constraints is fully captured by $G_{\mathrm{par}}^V$. Thus, $f_V(y)$ satisfies precisely the independences that derive from (1.2) and no others.

Methods have been developed to decide for any nonempty sets $\alpha$ and $\beta$ whether a given parent graph implies that $\alpha \perp\!\!\!\perp \beta|C$ holds for all distributions generated over it. Such methods have been called *separation criteria* because they check if the conditioning set $C$ is a separating set, in the sense of graph theory. Two quite different but equivalent criteria have been derived. Both are based on special types of paths in independence graphs. The first, by Geiger, Verma and Pearl [5] has been called d-separation, because it is applied to a directed acyclic graph. The second, by Lauritzen et al. [6], uses the basic graph theoretic notion of separation in an undirected graph. Such a graph, derived from $G_{\mathrm{par}}^V$, has been named a moral graph.

In this paper we use a different approach by associating first a joint Gaussian distribution of $Y$ with $G_{\mathrm{par}}^V$. For this distribution, the list of independencies (1.2) is equivalent to a set of zero population least-squares regression coefficients, that is, to a set of linear independencies

$$i \perp\!\!\!\perp j|\,\mathrm{par}_i \iff \beta_{i|j.\,\mathrm{par}_i} = 0,$$

where $\beta_{i|j.\,\mathrm{par}_i}$ is an adaption of the Yule–Cochran notation for the regression coefficient of $Y_j$, here in linear least-squares regression of $Y_i$ on $Y_{\mathrm{par}_i}$ and $Y_j$ for $j > i$ not a parent of $i$.

There are then two key results: (a) in distributions of arbitrary form generated over $G^V_{\mathrm{par}}$, probabilistic independence statements combine in the same way as linear independence statements; and (b) special types of path in $G^V_{\mathrm{par}}$ lead to dependence of $Y_\alpha$ on $Y_\beta$, given $Y_C$ in the *relevant family of Gaussian distributions generated over a given parent graph*, that is, in Gaussian distributions with parameters constrained only by the defining list of independencies (1.2) and having nonvanishing dependence of $Y_i$ on $Y_j$, given $Y_{\mathrm{par}\setminus j}$ for every $i \prec\!\!-j$ arrow present in the parent graph.

In any Gaussian distribution of $Y$, $\alpha \perp\!\!\!\perp \beta | C$ holds if and only if the population *coefficient matrix* of $Y_\beta$ in linear least-squares regression of $Y_\alpha$ on both $Y_\beta$ and $Y_C$ is zero. This matrix is related to linear equations associated with $G^V_{\mathrm{par}}$ using a generalization of the sweep operator for symmetric matrices [4] called *partial inversion*.

With another operator for binary matrices, named *partial closure*, so-called structural zeros in this matrix are expressed in terms of a special binary matrix representation derived from the parent graph. A particular zero submatrix will be shown to imply that $\alpha \perp\!\!\!\perp \beta | C$ holds in all distributions generated over $G^V_{\mathrm{par}}$.

This matrix criterion leads to a further path-based criterion for separation in directed acyclic graphs. Finally, equivalence of the new criterion to each of the two known separation criteria is established after having given first equivalent matrix formulations to each of these latter two path-based criteria.

**2. Edge matrices and induced independence statements.** Every independence graph has a matrix representation called its *edge matrix* (see [10, 12]). In this paper we are concerned with edge matrices derived from that of the parent graph called *induced edge matrices*.

2.1. *Edge matrix of the parent graph and linear recursive regressions.* The edge matrix of a parent graph is a $d \times d$ upper triangular binary matrix $\mathcal{A} = (\mathcal{A}_{ij})$ such that

$$(2.1) \qquad \mathcal{A}_{ij} = \begin{cases} 1, & \text{if and only of } i \prec\!\!-j \text{ in } G^V_{\mathrm{par}} \text{ or } i = j, \\ 0, & \text{otherwise.} \end{cases}$$

This matrix is the transpose of the usual adjacency matrix representation of $G^V_{\mathrm{par}}$, with additional ones along the diagonal. Because of the triangular form of the edge matrix $\mathcal{A}$, densities (1.1) generated over $G^V_{\mathrm{par}}$ are called *triangular systems*.

If the mean-centred random vector $Y$ generated over the parent graph has a joint Gaussian distribution, then each factor $f_{i|\mathrm{par}_i}(y_i|y_{\mathrm{par}_i})$ of equation (1.1) is a linear least-squares regression

$$Y_i = \sum_{j \in \mathrm{par}_i} \beta_{i|j.\,\mathrm{par}_i \setminus j} Y_j + \varepsilon_i,$$

where the residuals $\varepsilon_i$ are mutually independent with zero mean and variance $\sigma_{ii|\mathrm{par}_i}$. Then, the joint model can be written in the form

$$(2.2) \qquad\qquad\qquad AY = \varepsilon,$$

where $A$ is a real-valued $d \times d$ upper triangular matrix with ones along the diagonal. The $d \times 1$ vector $\varepsilon$ has zero mean and $\mathrm{cov}(\varepsilon) = \Delta$, a diagonal matrix with elements $\Delta_{ii} = \sigma_{ii|\mathrm{par}_i} > 0$. The covariance and concentration matrix of $Y$ are then, respectively, $\Sigma = A^{-1}\Delta A^{-\mathrm{T}}$ and $\Sigma^{-1} = A^{\mathrm{T}}\Delta^{-1}A$. An upper off-diagonal element of the generating matrix $A$ is $A_{ij} = -\beta_{i|j.\,\mathrm{par}_i \setminus j} \neq 0$ if the $ij$-arrow is present in the parent graph and $A_{ij} = 0$ if the $ij$-arrow is missing in $G^V_{\mathrm{par}}$.

In this Gaussian model, there is a one-to-one correspondence between missing $ij$-edges in the parent graph and zero parameters $A_{ij} = 0$. As a consequence, any such zero coincides in this case with a *structural zero*, that is a zero that holds for the relevant family of Gaussian distributions generated over $G^V_{\mathrm{par}}$. Therefore, the edge matrix $\mathcal{A}$ of the parent graph can be interpreted as the indicator matrix of zeros in $A$, that is $\mathcal{A} = \mathrm{In}[A]$, where the $\mathrm{In}[\cdot]$ operator transforms every nonzero entry of a matrix to be equal to one.

The edge matrix $\mathcal{A}$ of $G^V_{\mathrm{par}}$ is the starting point to find induced conditional independencies satisfied by all distributions generated over the parent graph. As we shall see, for any given Gaussian distribution generated over $G^V_{\mathrm{par}}$, independence statements are reflected by zeros in a matrix derived from $A$. Some of these zeros may be due to specific parametric constellations, others are consequences of the defining list of independencies (1.2). These latter zeros, that is the structural zeros in this matrix, show as zeros in edge matrices derived from $\mathcal{A}$, which in turn are matrix representations of *induced independence graphs*.

2.2. *Independence and structural zeros.* For any partitioning of the vertex set $V$ into node sets $M$, $\alpha$, $\beta$, $C$, where only $M$ and $C$ may be empty, a joint density that factorizes according to (1.1) is now considered in the form

$$f_V = f_{M|\alpha\beta C} f_{\alpha|\beta C} f_{\beta|C} f_C.$$

Marginalizing over $M$, as well as conditioning on $C$ and removing the corresponding variables, gives $f_{\alpha\beta|C} = f_{\alpha|\beta C} f_{\beta|C}$, so that $\alpha \perp\!\!\!\perp \beta | C$ holds if and only if $f_{\alpha|\beta C} = f_{\alpha|C}$.

For a Gaussian distribution, the independence $\alpha \perp\!\!\!\perp \beta | C$ is equivalently captured by

$$\Pi_{\alpha|\beta.C} = 0 \Longleftrightarrow \Sigma_{\alpha\beta|C} = 0,$$

where $\Pi_{\alpha|\beta.C}$ denotes the coefficient matrix of $Y_\beta$ in linear least-squares regression of $Y_\alpha$ on $Y_\beta$ and $Y_C$, and $\Sigma_{\alpha\beta|C}$ the joint conditional covariance matrix of $Y_\alpha$ and $Y_\beta$ given $Y_C$ (see Appendix A.2). The essential point is then that independence of $Y_\alpha$ and $Y_\beta$, given $Y_C$, is implied for all distributions generated over $G_{\text{par}}^V$, if and only if for Gaussian models (2.2) the induced coefficient matrix $\Pi_{\alpha|\beta.C}$ is implied to vanish for all permissible values of the parameters.

With $\alpha = a \setminus M$ and $\beta = b \setminus C$, we have $(\Pi_{a|b})_{\alpha,\beta} = \Pi_{\alpha|\beta.C}$. Therefore, the specific form of $\Pi_{\alpha|\beta.C}$ implied by $A$ in equation (2.2) can be derived by linear least-squares regression of $Y_a$ on $Y_b$ and the independence structure implied for it by $\mathcal{A}$ in equation (2.1) in a related way. Before turning to these tasks, we summarize how linear independence statements relate to probabilistic independencies specified with a parent graph.

2.3. *Combining independencies in triangular systems of densities.* It has been noted by Smith ([8], Example 2.8) that probabilistic and linear independencies combine in the same way. We prove a similar property, using assertions that have been discussed recently by Studený [9], where we take $V$ to be partitioned into $a$, $b$, $c$, $d$.

LEMMA 1. *In densities of arbitrary form generated over $G_{\text{par}}^V$, conditional independence statements combine as in a nondegenerate Gaussian distribution. This means that they satisfy the following statements, where we write, for instance, bc for the union of b and c:*

   (i) *symmetry: $a \perp\!\!\!\perp b | c$ implies $b \perp\!\!\!\perp a | c$;*
   (ii) *decomposition: $a \perp\!\!\!\perp bc | d$ implies $a \perp\!\!\!\perp b | d$;*
   (iii) *weak union: $a \perp\!\!\!\perp bc | d$ implies $a \perp\!\!\!\perp b | cd$;*
   (iv) *contraction: $a \perp\!\!\!\perp b | c$ and $a \perp\!\!\!\perp d | bc$ imply $a \perp\!\!\!\perp bd | c$;*
   (v) *intersection: $a \perp\!\!\!\perp b | cd$ and $a \perp\!\!\!\perp c | bd$ imply $a \perp\!\!\!\perp bc | d$;*
   (vi) *composition: $a \perp\!\!\!\perp c | d$ and $b \perp\!\!\!\perp c | d$ imply $ab \perp\!\!\!\perp c | d$.*

PROOF. The first four statements are basic properties of probabilities (see, e.g. [2]). Densities (1.1) generated over a parent graph also satisfy properties (v) and (vi), due to the full ordering of the node set the $ij$-dependence if and only if there is an $ij$-arrow in $G_{\text{par}}^V$ and the lack of any other constraint on the density.

For (v), the generating process implies for two nodes $i < j$ that of $a \perp\!\!\!\perp i | jd$ and $a \perp\!\!\!\perp j | id$, the statement $a \perp\!\!\!\perp j | id$ is not in the defining list (1.2) unless there are additional independencies. If $a \perp\!\!\!\perp j|id$ is to be satisfied, then at least $a \perp\!\!\!\perp j|d$ has to be in the defining list as well. And, in this case, $f_{ija|d} = f_{i|jd}f_{j|d} = f_{ij|d}$, so that $a \perp\!\!\!\perp ij|d$ is implied.

For (vi) and again $i < j$, both of $i \perp\!\!\!\perp c|d$ and $j \perp\!\!\!\perp c|d$ can only be in the defining list of independencies if the statement $i \perp\!\!\!\perp c|jd$ is also satisfied. And, in this case, $f_{ijc|d} = f_{i|jd}f_{j|d} = f_{ij|d}$, so that $ij \perp\!\!\!\perp c|d$ is implied. Equivalence of each of the assertions to statements involving least-squares coefficient matrices, proved in Lemma A.1, Appendix A.2, completes the proof. $\square$

2.4. *Partially inverted concentration matrices.* The induced parameter matrix $\Pi_{\alpha|\beta.C}$ is to be expressed in terms of the original parametrization $(A, \Delta)$. This is achieved in terms of the matrix operator partial inversion (see Appendix A.1 for a detailed summary of some of its properties).

Partial inversion with respect to any rows and columns $a$ of the concentration matrix, partitioned as $(a, b)$, transforms $\Sigma^{-1}$ into $\mathrm{inv}_a \Sigma^{-1}$

$$(2.3) \qquad \Sigma^{-1} = \begin{pmatrix} \Sigma^{aa} & \Sigma^{ab} \\ \cdot & \Sigma^{bb} \end{pmatrix}, \qquad \mathrm{inv}_a \Sigma^{-1} = \begin{pmatrix} \Sigma_{aa|b} & \Pi_{a|b} \\ \sim & \Sigma^{bb.a} \end{pmatrix},$$

where the $\cdot$ notation indicates entries in a symmetric matrix, and the $\sim$ notation denotes entries in a matrix that is symmetric except for the sign. The submatrix $\Pi_{a|b}$ is as defined before; submatrix $\Sigma_{aa|b} = (\Sigma^{aa})^{-1}$ is the covariance matrix of $Y_a - \Pi_{a|b}Y_b$ and submatrix $\Sigma^{bb.a} = \Sigma_{bb}^{-1}$ is the marginal concentration matrix of $Y_b$. We denote by $\tilde{A}$ the accordingly partitioned matrix $A$ in which the original order is preserved both within $a$ and within $b$, but which is typically asymmetric and not triangular.

An important property of partial inversion is that

$$\tilde{A}\begin{pmatrix} Y_a \\ Y_b \end{pmatrix} = \begin{pmatrix} \varepsilon_a \\ \varepsilon_b \end{pmatrix} \quad \text{implies} \quad \mathrm{inv}_a \tilde{A}\begin{pmatrix} \varepsilon_a \\ Y_b \end{pmatrix} = \begin{pmatrix} Y_a \\ \varepsilon_b \end{pmatrix},$$

so that, with $B = \mathrm{inv}_a \tilde{A}$, one obtains directly the equations in $Y_b$ from which $Y_a$ has been removed as

$$(2.4) \qquad B_{bb}Y_b = \eta_b, \qquad \eta_b = \varepsilon_b - B_{ba}\varepsilon_a.$$

In addition, direct covariance computations give

$$\tau = \mathrm{cov}\begin{pmatrix} \varepsilon_a \\ \eta_b \end{pmatrix} = \begin{pmatrix} \Delta_{aa} & -\Delta_{aa}B_{ba}^{\mathrm{T}} \\ \cdot & \Delta_{bb} + B_{ba}\Delta_{aa}B_{ba}^{\mathrm{T}} \end{pmatrix}.$$

Therefore, equations in $Y_a$ corrected for linear dependence on $Y_b$ and having residuals uncorrelated with $\eta_b$ are, with $H = \mathrm{inv}_b \tau$,

$$(2.5) \qquad Y_a - \Pi_{a|b}Y_b = B_{aa}\eta_a, \qquad \eta_a = \varepsilon_a - H_{ab}\eta_b.$$

Lemma 2 below is now a direct consequence of the definition of $H$ and equations (2.4) and (2.5). It leads, after expansion of the matrix $H$, to an explicit expression of the matrix of least-squares regression coefficients $\Pi_{a|b}$ as a function of $\Delta$ and $B$, used later for Proposition 2. Similarly, the explicit expression of $\Sigma_{aa|b}$ will be used for Proposition 3 and of $\Sigma^{bb.a}$ for Proposition 4.

LEMMA 2 (Wermuth and Cox [10]). *For a linear triangular system (2.2), with* $\Sigma = \text{cov}(Y)$, *a any subset of* $V$, $b = V \setminus a$ *and* $H$ *as defined for equations (2.4), (2.5),*

$$(2.6) \qquad \text{inv}_a \, \Sigma^{-1} = \begin{pmatrix} B_{aa} H_{aa} B_{aa}^{\mathrm{T}} & B_{ab} + B_{aa} H_{ab} B_{bb} \\ \sim & B_{bb}^{\mathrm{T}} H_{bb} B_{bb} \end{pmatrix}.$$

2.5. *Induced edge matrices.* To obtain induced edge matrices, we display first explicitly $B = \text{inv}_a \, \tilde{A}$, and $\mathcal{B} = \text{zer}_a \, \tilde{A}$, obtained by what is called the operator of partial closure (see Appendix A.1 for some of its properties). It finds both the structural zeros in $B$ and the edge matrix induced by $\mathcal{A}$ for what we define below as the partial ancestor graph, with respect to subset $a$ of node set $V$.

$$(2.7) \qquad \begin{aligned} B &= \begin{pmatrix} A_{aa}^{-1} & -A_{aa}^{-1}\tilde{A}_{ab} \\ \tilde{A}_{ba}A_{aa}^{-1} & A_{bb} - \tilde{A}_{ba}A_{aa}^{-1}\tilde{A}_{ab} \end{pmatrix}, \\ \mathcal{B} &= \text{In}\left[ \begin{pmatrix} \mathcal{A}_{aa}^- & \mathcal{A}_{aa}^-\tilde{\mathcal{A}}_{ab} \\ \tilde{\mathcal{A}}_{ba}\mathcal{A}_{aa}^- & \mathcal{A}_{bb} + \tilde{\mathcal{A}}_{ba}\mathcal{A}_{aa}^-\tilde{\mathcal{A}}_{ab} \end{pmatrix} \right], \end{aligned}$$

with

$$\mathcal{A}_{aa}^- = \text{In}[(k\mathcal{I}_{aa} - \mathcal{A}_{aa})^{-1}],$$

where $\mathcal{I}_{aa}$ denotes an identity matrix of dimension $d_a$ and $k = d_a + 1$. The matrix $\mathcal{A}_{aa}^-$ provides the structural zeros in $A_{aa}^{-1}$ and the edge matrix of the transitive closure of the graph with edge matrix $\mathcal{A}_{aa}$, for which fast algorithms are also available (see [3]).

The transition from a matrix of parameters in a linear system, such as $B$ in equation (2.7), to a corresponding induced edge matrix, $\mathcal{B}$, is generalized with Lemma 3 below.

LEMMA 3 (Wermuth and Cox [11]). *Let induced parameter matrices be defined by parameter components of a linear system of the type* $FY = \zeta$ *with possibly correlated residuals* $\zeta$, *so that the defining matrix products hide no self-cancellation of an operation such as a matrix multiplied by its inverse. Further let the structural zeros of* $F$ *be given by* $\mathcal{F}$. *Then, the induced edge matrix components are obtained by replacing, in a given sum of products:*

(i) *every inverse matrix, say* $F_{aa}^{-1}$ *by the binary matrix of its structural zeros* $\mathcal{F}_{aa}^-$;

(ii) *every diagonal matrix by an identity matrix of the same dimension;*
(iii) *every other submatrix, say* $-F_{ab}$ *or* $F_{ab}$, *by the corresponding binary submatrix of structural zeros,* $\mathcal{F}_{ab}$;

*and then applying the indicator function.*

By using Lemma 3, each submatrix of a linear parameter matrix is substituted by a nonnegative matrix having the appropriate structural zeros. By multiplying, summing and applying the indicator function, all structural zeros are preserved and no additional zeros are generated, while some structural zeros present in $\mathcal{F}$ may be changed to ones.

After applying Lemma 3 to equation (2.6), the edge matrix components induced by $\mathcal{A}$ for $\mathrm{inv}_a \Sigma^{-1}$ result. These components are $\mathcal{P}_{a|b}$ of $\Pi_{a|b}$, $\mathcal{S}_{aa|b}$ of $\Sigma_{aa|b}$ and $\mathcal{S}^{bb.a}$ of $\Sigma^{bb.a}$.

LEMMA 4 (Wermuth, Wiedenbeck and Cox [12]).  *The edge matrix components induced by a parent graph for* $\mathrm{inv}_a \Sigma^{-1}$ *in (2.3) are*

$$(2.8) \qquad \begin{pmatrix} \mathcal{S}_{aa|b} & \mathcal{P}_{a|b} \\ . & \mathcal{S}^{bb.a} \end{pmatrix} = \mathrm{In}\left[ \begin{pmatrix} \mathcal{B}_{aa}\mathcal{H}_{aa}\mathcal{B}_{aa}^{\mathrm{T}} & \mathcal{B}_{ab} + \mathcal{B}_{aa}\mathcal{H}_{ab}\mathcal{B}_{bb} \\ . & \mathcal{B}_{bb}^{\mathrm{T}}\mathcal{H}_{bb}\mathcal{B}_{bb} \end{pmatrix} \right],$$

*where*

$$\mathcal{H} = \mathrm{zer}_b \begin{pmatrix} \mathcal{I}_{aa} & \mathcal{B}_{ba}^{\mathrm{T}} \\ . & \mathcal{I}_{bb} + \mathcal{B}_{ba}\mathcal{B}_{ba}^{\mathrm{T}} \end{pmatrix}.$$

Lemma 4 leads to the following matrix criteria for independencies implied by a parent graph, where $V$ is partitioned as before into $M$, $\alpha$, $\beta$, $C$.

PROPOSITION 1.  *The parent graph* $G_{\mathrm{par}}^V$ *implies, for every density generated over it that:*

(i) $\alpha \perp\!\!\!\perp \beta | C$ *holds if and only if* $(\mathcal{P}_{a|b})_{\alpha,\beta} = \mathcal{P}_{\alpha|\beta.C} = 0$;
(ii) $\alpha \perp\!\!\!\perp M | b$ *holds if and only if* $\mathcal{S}_{\alpha M|b} = 0$;
(iii) $\beta \perp\!\!\!\perp C | a$ *holds if and only if* $\mathcal{S}^{\beta C.a} = 0$.

PROOF.  For (i), if there is a one for $i, j$ in $\mathcal{P}_{\alpha|\beta.C}$, then $\beta_{i|j.C\beta\setminus j} \neq 0$ holds in the relevant family of Gaussian densities generated over $G_{\mathrm{par}}^V$ (see, e.g. [11]). If, instead, the $ij$-edge is missing, then $i \perp\!\!\!\perp j | C\beta$ is implied for every member of the family of Gaussian distributions. For $\mathcal{P}_{\alpha|\beta.C} = 0$, the statement $\alpha \perp\!\!\!\perp \beta | C$ results with Lemma 1 for every distribution generated over $G_{\mathrm{par}}^V$. For (ii) and (iii), the arguments are analogous.  □

**3. A path-based interpretation of the matrix criterion.** To give a path interpretation of the matrix criterion $\mathcal{P}_{\alpha|\beta.C} = 0$, we first summarize some definitions related to paths and graphs.

Two nodes $i$ and $j$ in an independence graph have at most one edge. If the $ij$-edge is present in the graph, then the node pair $i, j$ is *coupled*; if the $ij$-edge is missing, the node pair is *uncoupled*. An $ij$-*path* connects the path endpoints $i$ and $j$ by a sequence of edges visiting distinct nodes. All nodes of a path except for the endpoint nodes are called the *inner nodes* of the path. An edge is regarded as a path without inner nodes. For a graph in node set $V$ and $a \subset V$, the *subgraph induced by $a$* is obtained by removing all nodes and edges outside $a$.

Both a graph and a path are called *directed* if all its edges are arrows. Directed graphs can have the following V-configurations, that is subgraphs induced by three nodes and having two edges,

$$i \prec\!\!\!— t \prec\!\!\!— j, \qquad i \prec\!\!\!— s —\!\!\!\succ j, \qquad i —\!\!\!\succ c \prec\!\!\!— j,$$

where the inner node is called a transition ($t$), a source ($s$) and a collision node ($c$), respectively. A directed path is *direction-preserving* if all its inner nodes are transition nodes. If in a direction-preserving path an arrow starts at node $j$ and points to $i$, then node $j$ is an *ancestor* of $i$, node $i$ a *descendant* of $j$, and the $ij$-path is called a *descendant–ancestor path*.

Node $j$ is an *$a$-line ancestor* of node $i$ if all inner nodes in the descendant–ancestor $ij$-path are in set $a$. A directed path is an *alternating path* if it has at least one inner node and the direction of the arrows changes at each inner node. This implies that no inner node is a transition node and that the inner nodes alternate between source and collision nodes. A parent graph is said to be *transitive* if it contains no transition-oriented V-configuration or, equivalently, if for each node $i$ the set $\mathrm{par}_i$ of its parents coincides with its set of ancestors.

The *partial ancestor graph*, with respect to nodes $a$, denoted by $G_{\mathrm{anc}}^{V.a}$, is an induced graph defined, for a reordered node set, by the edge matrix $\mathcal{B} = \mathrm{zer}_a \tilde{\mathcal{A}}$ of equation (2.7). The elements of $\mathcal{B}$ are equivalently given by

$$\mathcal{B}_{ij} = \begin{cases} 1, & \text{if and only if } j \text{ is an } a\text{-line ancestor of } i \text{ in } G_{\mathrm{par}}^V \text{ or } i = j, \\ 0, & \text{otherwise,} \end{cases}$$

(3.1)

with nodes ordered as for equation (2.7). Since $\mathcal{B}$ in (3.1) implies that every $a$-line descendant–ancestor path in $G_{\mathrm{par}}^V$ is in $G_{\mathrm{anc}}^V$ closed by an arrow that points to the descendant, the corresponding operator has been named partial closure. Induced edge matrices and induced linear parameter matrices may be calculated within the statistical environment R (see [7]).

Figure 1(a) shows a parent graph in six nodes, Figure 1(b) its partial ancestor graph with respect to $a = \{1, 2, 3\}$, and Figure 1(c) the induced graph for the conditional dependence of $Y_a$ given $Y_b$. In this example, $(a, b)$ is an *order compatible partition of the node set $V$*, that is a mere split of $V = (1, \ldots, d)$ into two components without any reordering. For such a split, $\mathcal{A}_{ba} = \mathcal{B}_{ba} = 0$ implies $\mathcal{P}_{a|b} = \mathcal{B}_{ab}$. With a transitive parent graph, in addition, $\mathcal{P}_{a|b} = \mathcal{A}_{ab}$ for all order compatible partitionings of $V$.

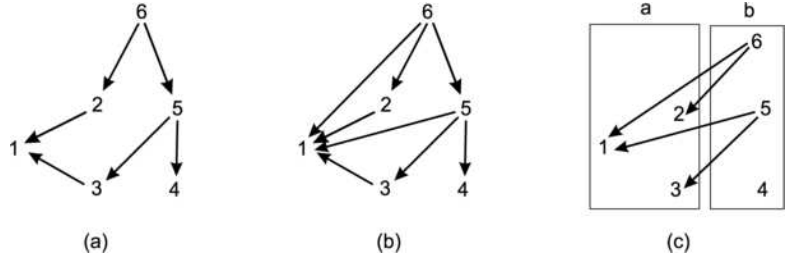

FIG. 1.  (a) *A parent graph $G_{\mathrm{par}}^V$ in six nodes.* (b) *Its partial ancestor graph $G_{\mathrm{anc}}^{V.a}$ with $a = \{1,2,3\}$.* (c) *The induced graph with edge matrix $\mathcal{P}_{a|b}$ for the conditional dependence of $Y_a$ on $Y_b$, where $a = \{1,2,3\}$ and $b = \{4,5,6\}$.*

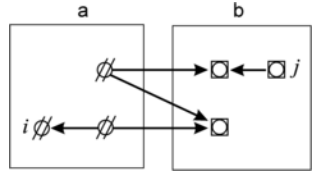

FIG. 2.  *An alternating path in $G_{\mathrm{anc}}^{V.a}$ from b to a; active since, of its inner nodes, each source node is in* a *and each collision node in* $b = V \setminus a$*; nodes in* a *are indicated by* ∅, *those in* b *by* ▢.

If, instead, $a$ is an arbitrary subset of $V$, then the graph with edge matrix $\mathcal{P}_{a|b}$ may contain additional edges compared to the graph with edge matrix $\mathcal{B}_{ab}$ [see equation (2.8)]. Such edges are due to the following type of path.

DEFINITION 3.1.   An $ij$-alternating path in the partial ancestor graph $G_{\mathrm{anc}}^{V.a}$ is called active if of its inner nodes every collision node is in $b$ and every source node is in $a$.

Thus, every off-diagonal one in the edge matrix (2.8) induced by a parent graph for $\mathrm{inv}_a \Sigma^{-1}$ can be identified in the partial ancestor graph by what we call its active paths. In diagrams of paths, we indicate nodes within $a$ as crossed out, ∅, and nodes within $b$ as boxed in, ▢, such as in Figure 2.

DEFINITION 3.2.   An $ij$-path in the partial ancestor graph $G_{\mathrm{anc}}^{V.a}$ is active if it is an $ij$-edge or it is an active alternating path.

PROPOSITION 2. *For node set $V$, partitioned into $a$ and $b$ and having node $i$ in $a$ and node $j$ in $b$, the induced graph with edge matrix $\mathcal{P}_{a|b}$ has an $ij$-arrow if and only if there is an active $ij$-path in the partial ancestor graph, with respect to $a$.*

The essence of the proof is the expansion of the sum of products defining $\mathcal{P}_{a|b}$ in equation (2.8) into submatrices of $\mathcal{B} = \mathrm{zer}_a \tilde{\mathcal{A}}$ and the interpretation of each matrix operation in terms of arrows present in $G_{\mathrm{anc}}^{V.a}$.

PROOF OF PROPOSITION 2. From equation (2.8) defining $\mathcal{P}_{a|b}$, there is an $ij$-one in $\mathcal{P}_{a|b}$ if and only if

$$\mathcal{B}_{ij} = 1 \quad \text{or} \quad \mathcal{B}_{ia}\mathcal{B}_{ba}^{\mathrm{T}}(\mathcal{I}_{bb} + \mathcal{B}_{ba}\mathcal{B}_{ba}^{\mathrm{T}})^{-}\mathcal{B}_{bj} \geq 1.$$

The first condition $\mathcal{B}_{ij} = 1$ holds if there is an arrow pointing from $j$ in $b$ to $i$ in $a$. The second condition holds if either an arrow points from $a$ to $b$, or if $i$ and $j$ are uncoupled but connected by an active alternating path. This interpretation of the second condition as an active alternating path is illustrated with the following scheme, in which the inner nodes are shown by their location in either $a$ or $b$:

$$
\begin{array}{ccccc}
\text{Edge matrix} & \mathcal{B}_{ia} & \mathcal{B}_{ba}^{\mathrm{T}} & (\mathcal{I}_{bb} + \mathcal{B}_{ba}\mathcal{B}_{ba}^{\mathrm{T}})^{-} & \mathcal{B}_{bj} \\
\text{Path} & {}_{i \leftarrow a} & {}_{a \rightarrow b} \ {}_{b} \searrow & {}_{a} \nearrow {}^{b} \dots {}^{b} \searrow & {}_{a} \nearrow {}^{b} \ {}_{b \leftarrow j}
\end{array}.
$$

Such an $ij$-path induces a dependence of $Y_i$ on $Y_j$ given $Y_{b \setminus j}$ in the relevant family of Gaussian distributions generated over $G_{\mathrm{par}}^{V}$.  □

The scheme shows that for the second condition, there has to be at least one arrow in $G_{\mathrm{anc}}^{V.a}$ pointing from $a$ to $b$. Therefore, for any order compatible split of $V$ into $(a, b)$ there is no active alternating path in $G_{\mathrm{anc}}^{V.a}$.

Proposition 2 leads to the following path criterion.

CRITERION 1. *If there is no active path between $\alpha$ and $\beta$ in the partial ancestor graph $G_{\mathrm{anc}}^{V.a}$, then $\alpha \perp\!\!\!\perp \beta | C$ in every joint density generated over the given parent graph.*

Figure 3(a) shows a parent graph and illustrates the use of this criterion.

**4. Equivalence to known separation criteria.** For a discussion of the criteria available in the literature to verify whether $\alpha \perp\!\!\!\perp \beta | C$ is implied by a given parent graph $G_{\mathrm{par}}^{V}$, we take throughout the node set $V$ to be partitioned into $M, \alpha, \beta, C$, where only $M$ and $C$ may be empty, and every node pair between $\alpha$ and $\beta$ to be uncoupled in $G_{\mathrm{par}}^{V}$.

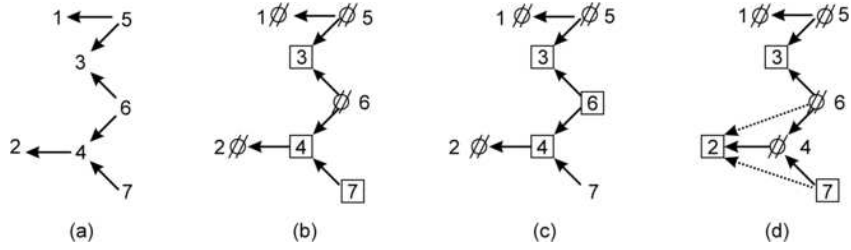

FIG. 3. *Illustration of Criterion 1.* (a) *Parent graph.* (b)–(d) *Is* $\alpha \perp\!\!\!\perp \beta | C$ *implied for* $\alpha = \{5\}$, $\beta = \{7\}$ *and different choices of* $C$; $\{\boxed{\phantom{o}}\} = \beta \cup C$ *and* $\{\emptyset\!\!\!/\} = \alpha \cup M$. (b) $5 \perp\!\!\!\perp 7 | \{3, 4\}$ *not implied, since the alternating path* $(5, 3, 6, 4, 7)$ *in* $G_{\mathrm{anc}}^{V,a}$ *is active with its inner source node,* $6$, *in* $\{\emptyset\!\!\!/\}$ *and its inner collision nodes,* $3, 4$ *in* $\{\boxed{\phantom{o}}\}$; (c) $5 \perp\!\!\!\perp 7 | \{3, 4, 6\}$ *implied, since inner source node* $6$ *in* $\{\boxed{\phantom{o}}\}$; (d) $5 \perp\!\!\!\perp 7 | \{2, 3\}$ *not implied, since the alternating path* $(5, 3, 6, 2, 7)$ *in* $G_{\mathrm{anc}}^{V,a}$ *is active.*

One known criterion uses definitions of a d-connecting path and of d-separation, where the letter d is to remind us that the definitions pertain to a directed acyclic graph. In $G_{\mathrm{par}}^V$, a path is said to be *d-connecting relative to* $C$ if along it every inner collision node is in $C$ or has a descendant in $C$ and every other inner node is outside $C$. And, two disjoint sets of nodes $\alpha$ and $\beta$ are said to be *d-separated* by $C$ if and only if, relative to $C$, there is no d-connecting *path between* $\alpha$ *and* $\beta$ that is between a node in $\alpha$ and a node in $\beta$.

CRITERION 2 (Geiger, Verma and Pearl [5]). *If* $\alpha$ *and* $\beta$ *are d-separated by* $C$ *in the parent graph, then* $\alpha \perp\!\!\!\perp \beta | C$ *in every joint density generated over the given parent graph.*

Proposition 3 below gives a matrix criterion that we will show to be equivalent to d-separation. For this, we denote by $g = V \setminus C$ the union of $\alpha, \beta$ and $M$, and by $\mathcal{F} = \mathrm{zer}_g \tilde{A}$ the edge matrix of $G_{\mathrm{anc}}^{V,g}$, the ancestor graph with respect to $g$. By equating $a$ to $g$ and $b$ to $C$, equation (2.8) gives

$$\mathcal{S}_{gg|C} = \mathrm{In}[\mathcal{F}_{gg}(\mathcal{I}_{gg} + \mathcal{F}_{Cg}^{\mathrm{T}} \mathcal{F}_{Cg})^{-} \mathcal{F}_{gg}^{\mathrm{T}}]$$

as the edge matrix induced by $\mathcal{A}$ for $\Sigma_{gg|C}$. Edges in the corresponding undirected graph, named the induced *covariance graph* of $Y_g$ given $Y_C$ [10], are drawn later in Figure 5 as dashed lines.

PROPOSITION 3. *In the parent graph* $G_{\mathrm{par}}^V$, *sets* $\alpha$ *and* $\beta$ *are d-separated by* $C$ *if and only if* $\mathcal{S}_{\alpha\beta|C} = 0$, *where* $\mathcal{S}_{\alpha\beta|C}$ *is the edge matrix induced by* $\mathcal{A}$ *for* $\Sigma_{\alpha\beta|C}$.

The following scheme shows that $\mathcal{S}_{\alpha\beta|C} = 0$ means the absence of any active path in $G_{\mathrm{anc}}^{V.g}$ from $\alpha$ to $\beta$, both subsets of $g$:

Edge matrix $\quad \mathcal{F}_{\alpha g} \quad (\mathcal{I}_{gg} + \mathcal{F}_{Cg}^{\mathrm{T}}\mathcal{F}_{Cg})^{-} \quad \mathcal{F}_{\beta g}^{\mathrm{T}}$

Path $\qquad {}^{\alpha \leftarrow g}\,{}^{g}\searrow_{C}\swarrow^{g}\,\ldots\,{}^{g}\searrow_{C}\swarrow^{g}\,{}^{g \to \beta}.$

PROOF OF PROPOSITION 3. By the definition of partial inversion and of $\mathcal{S}_{gg|C}$, a d-connecting $ij$-path in $G_{\mathrm{par}}^{V}$ relative to $C$ and without inner collision nodes, forms in $G_{\mathrm{anc}}^{V.g}$ an $ij$-edge or a sequence of three nodes $(i, s, j)$, with $s$ a source in $g$. Both types are active paths in $G_{\mathrm{anc}}^{V.g}$.

If there is a d-connecting $ij$-path in $G_{\mathrm{par}}^{V}$, relative to $C$ and having an inner collision node, then $\mathrm{zer}_g\,\bar{\mathcal{A}}$ generates an active alternating $ij$-path in $G_{\mathrm{anc}}^{V.g}$ as follows. Every inner source node and every inner collision node within $C$ is preserved. Every inner collision $h$ outside $C$ is replaced by its first $g$-line descendant $h_C$ within $C$. Every transition node $t$ in an inner node sequence $(i, t, j)$ is removed via the $ij$-edge present in $G_{\mathrm{anc}}^{V.g}$.

Conversely, if there is an active alternating $ij$-path in $G_{\mathrm{anc}}^{V.g}$, then, by these constructions, there is a d-connecting path relative to $C$ in $G_{\mathrm{par}}^{V}$. $\quad\square$

Figure 4 shows a d-connecting path relative to conditioning set $C = \{\boxdot\}$ and the corresponding active alternating path in $G_{\mathrm{anc}}^{V.g}$ with $g = \{\varnothing\}$.

The other criterion in the literature uses an undirected graph called the *moral graph* of $\alpha, \beta$, and $C$. This moral graph is constructed in three steps. One obtains the subgraph induced by the union of the nodes $\alpha, \beta$, and $C$ and their ancestors. One joins by a line every uncoupled pair of parents having a common offspring. One replaces every arrow in the resulting graph by a line. Then, the separation criterion for undirected graphs is used to give Criterion 3. In the moral graph, $C$ separates $\alpha$ from $\beta$ if every path from a node in $\alpha$ to one in $\beta$ has a node in $C$.

CRITERION 3 (Lauritzen, Dawid, Larsen and Leimer [6]). *If $\alpha$ and $\beta$ are separated by $C$ in the moral graph of $\alpha, \beta, C$, then $\alpha \perp\!\!\!\perp \beta | C$ in every joint density generated over the given parent graph.*

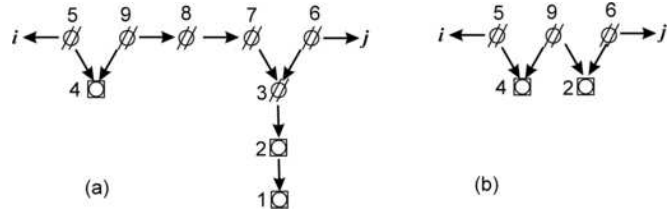

FIG. 4. (a) *Example of a d-connecting $ij$-path in $G_{\mathrm{par}}^{V}$ relative to $C$ with inner nodes* 5, 4, 9, 8, 7, 3, 6 *and* (b) *the corresponding active alternating $ij$-path in $G_{\mathrm{anc}}^{V.g}$ with inner nodes* 5, 4, 9, 2, 6.

By the definition of the moral graph and by equating $b$ to $Q$, the union of $\alpha, \beta$ and $C$ and their ancestors, and $a$ to $O = V \setminus Q$, equation (2.8) gives the edge matrix of the moral graph of $\alpha, \beta, C$ as $\mathcal{S}^{QQ.O}$ the edge matrix induced by $\mathcal{A}$ for $\Sigma^{QQ.O}$. Since there is no path leading from $O$ to $Q$, this edge matrix has the special form

$$(4.1) \qquad \mathcal{S}^{QQ.O} = \text{In}[\mathcal{A}_{QQ}^{\mathrm{T}} \mathcal{A}_{QQ}] \quad \text{due to} \quad \mathcal{A}_{OQ} = 0.$$

Thus, the induced graph contains an $ij$-edge if and only if in the parent graph either there is an $ij$-edge, or $\mathcal{A}_{hi}\mathcal{A}_{hj} = 1$ for some node $h < i < j$ in $Q$, that is for nodes $i$ and $j$ having a common offspring in $Q$.

Proposition 4 below gives a matrix criterion that we will show to be equivalent to separation in the moral graph of $\alpha, \beta$ and $C$. For this, we let $q = V \setminus M$, so that the set $q$ denotes the union of $\alpha, \beta$ and $C$. Furthermore, we denote by $\mathcal{Z} = \text{zer}_M \tilde{\mathcal{A}}$ the edge matrix of the induced partial ancestor graph with respect to $M$. Then, by equating $b$ to $q$ and $a$ to $M$, equation (2.8) gives, as edge matrix induced by $\mathcal{A}$ for $\Sigma^{qq.M}$,

$$(4.2) \qquad \mathcal{S}^{qq.M} = \text{In}[\mathcal{Z}_{qq}^{\mathrm{T}}(\mathcal{I}_{qq} + \mathcal{Z}_{qr}\mathcal{Z}_{qr}^{\mathrm{T}})^{-}\mathcal{Z}_{qq}],$$

where $r$ denotes the set of ancestors of $q$ within $M$. Again, the special form of $\mathcal{S}^{qq.M}$ is due to $\mathcal{A}_{OQ} = 0$. It leads to $\mathcal{Z}_{QQ} = (\text{zer}_r \tilde{\mathcal{A}})_{Q,Q}$ since $M = r \cup O$. As a consequence, also $\mathcal{S}^{qq.M} = \mathcal{S}^{qq.r}$. Edges in this type of undirected graph, named the *induced concentration graph* of $Y_q$, are drawn in Figure 5 as full lines.

PROPOSITION 4. *In the moral graph of $\alpha, \beta$ and $C$, set $\alpha$ is separated from $\beta$ by $C$ if and only if $\mathcal{S}^{\alpha\beta.M} = 0$, where $\mathcal{S}^{\alpha\beta.M}$ is the edge matrix induced by $\mathcal{A}$ for $\Sigma^{\alpha\beta.M}$.*

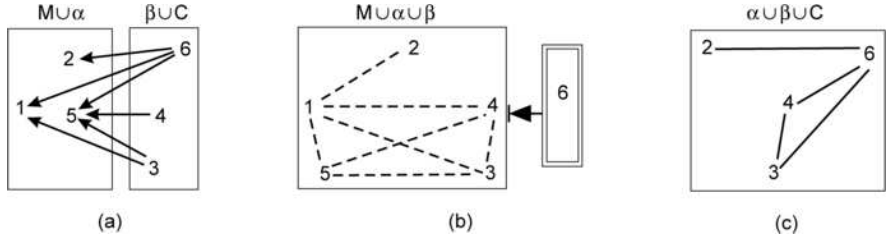

Fig. 5. *Graphs induced by the parent graph in Figure 1*(a)*, each of which shows that $\alpha \perp\!\!\!\perp \beta | C$ holds for $\alpha = \{2\}$, $\beta = \{3, 4\}$ and $C = \{6\}$ by all edges between $\alpha$ and $\beta$ being missing. The graph induced by $\mathcal{A}$* (a) *with edge matrix $\mathcal{P}_{a|b}$ for $a = V \setminus b$ and $b$ the union of $\beta$ and $C$;* (b) *with edge matrix $\mathcal{S}_{gg|C}$ for $g = V \setminus C$;* (c) *with edge matrix $\mathcal{S}^{qq.M}$ for $q = V \setminus M$.*

The following scheme shows that $\mathcal{S}^{\alpha\beta.M} = 0$ means the absence of any active path in $G^{V.q}_{\text{anc}}$ from $\alpha$ to $\beta$, which are both subsets of $q$:

Edge matrix $\qquad \mathcal{Z}^{\mathrm{T}}_{q\alpha} \qquad (\mathcal{I}_{qq} + \mathcal{Z}_{qr}\mathcal{Z}^{\mathrm{T}}_{qr})^{-} \qquad \mathcal{Z}_{q\beta}$

Path $\qquad\qquad {}_{\alpha\to q}\ {}_q \searrow\ {}_r \nearrow {}^q \dots {}^q \searrow\ {}_r \nearrow {}^q\ {}_{q\leftarrow\beta}$.

PROOF OF PROPOSITION 4. We show that the edge matrix $\mathcal{S}^{qq.M}$ can equivalently be obtained via edges and alternating paths in $G^{V.M}_{\text{anc}}$ and by closing all $r$-line paths in the moral graph for $\alpha, \beta, C$, which has the edge matrix $\mathcal{S}^{QQ.O}$ in equation (4.1).

For this, we note first that $Q = q \cup r$ and we recall that, for $\mathcal{S}^{QQ.O}$, all collision-oriented V-configurations in the parent graph are closed, that have a common collision node in $Q$. Then, in the resulting concentration graph, all $r$-line paths are closed by partial inversion of $\mathcal{S}^{QQ.O}$, with respect to $r$. This gives, for the subgraph induced by nodes $q$, the edge matrix $(\text{zer}_r \tilde{\mathcal{S}}^{QQ.O})_{q,q}$.

For the edge matrix $\mathcal{S}^{qq.M}$ in (4.2), all $r$-line ancestor–descendant paths are closed first with $\mathcal{Z} = \text{zer}_r \tilde{\mathcal{A}}$, whereby every collision node within $r$, each of which has a $q$-line descendant in $q$, is replaced by the first descendant in $q$ (see Figure 4 for an illustration). Then, the active alternating path in $G^{V.M}_{\text{anc}}$ has every source node in $r$ and every collision node in $q$.

Thus, $(\text{zer}_r \tilde{\mathcal{S}}^{QQ.O})_{q,q} = \mathcal{S}^{qq.M}$, since for both edge matrices exactly the following types of V-configurations are closed in the subgraph induced by $q$ in $G^V_{\text{par}}$:

$$i \leftarrow\!\!-\ r \prec\!\!-\!\!j, \qquad i \prec\!\!-\ r \to j \quad \text{and} \quad i \leftarrow\!\!-\ Q \to j. \qquad \square$$

Our final result establishes the equivalence of the three path-based separation criteria for $V$ partitioned as before into $M, \alpha, \beta, C$ and explains why proofs of equivalence become complex when they are based exclusively on paths that induce edges in different types of graph.

PROPOSITION 5. *The following assertions are equivalent. Between $\alpha$ and $\beta$ there is:*

(i) *an active path in the partial ancestor graph with respect to $M$ and $\alpha$;*

(ii) *a d-connecting path relative to $C$ in the parent graph;*

(iii) *a $M$-line path in the moral graph of $\alpha, \beta$ and $C$.*

PROOF. By using Propositions 2 to 4, the results follows after noting that

$$\mathcal{P}_{\alpha|\beta.C} = 0 \iff \mathcal{S}_{\alpha\beta|C} = 0 \iff \mathcal{S}^{\alpha\beta.M} = 0. \qquad \square$$

Figure 5 illustrates the result of Proposition 5 for the parent graph of Figure 1(a), with $\alpha = \{2\}$, $\beta = \{3,4\}$, $C = \{6\}$ and $M = \{1,5\}$. The independence statement $2 \perp\!\!\!\perp \{3,4\}|6$ is implied by the parent graph, since the subgraph induced by nodes $\alpha$ and $\beta$ has no edge between $\alpha$ and $\beta$ in the induced graph with edge matrix $\mathcal{P}_{a|b}$ in Figure 5(a), with edge matrix $\mathcal{S}_{gg|C}$ for $g = V \setminus C$ in Figure 5(b), and with edge matrix $\mathcal{S}^{qq.M}$ for $q = V \setminus M$ in Figure 5(c).

## APPENDIX: OPERATORS AND LINEAR INDEPENDENCIES

**A.1. Partial inversion and partial closure.** Two matrix operators, introduced and studied by Wermuth, Wiedenbeck and Cox [12], permit stepwise transformations of parameters in linear systems and edge matrices of independence graphs, respectively. Note that $M$ now denotes a matrix.

Let $M$ be a square matrix of dimension $d$, for which all principal submatrices are invertible. For a given integer $1 \leq k \leq d$, *partial inversion* of $M$ with respect to $k$, transforms $M$ into a matrix $N = \mathrm{inv}_k\, M$ of the same dimensions, where, for all $i, j \neq k$,

$$
\begin{aligned}
(A.1) \quad & N_{kk} = 1/M_{kk}, \\
& N_{ik} = M_{ik}/M_{kk}, \\
& N_{kj} = -M_{kj}/M_{kk}, \\
& N_{ij} = M_{ij} - M_{ik}M_{kj}/M_{kk}.
\end{aligned}
$$

Then the matrix $\mathcal{N}$ of structural zeros in $N = \mathrm{inv}_k\, M$ that remain after partial inversion of $M$ on $k$ is defined by $\mathcal{N} = \mathrm{zer}_k\, \mathcal{M}$:

$$
\begin{aligned}
(A.2) \quad & \mathcal{N}_{kk} = 1, \\
& \mathcal{N}_{ik} = \mathcal{M}_{ik}, \\
& \mathcal{N}_{kj} = \mathcal{M}_{kj}, \\
& \mathcal{N}_{ij} = \begin{cases} 1, & \text{if } \mathcal{M}_{ij} = 1 \text{ or } \mathcal{M}_{ik}\mathcal{M}_{kj} = 1, \\ 0, & \text{otherwise.} \end{cases}
\end{aligned}
$$

Partial inversion of $M$, with respect to a sequence of indices $a$, applies the operator (A.1) in sequence to all elements of $a$ and is denoted by $\mathrm{inv}_a\, M$. Similarly, partial closure (A.2) of $\mathcal{M}$, with respect to $a$, is denoted by $\mathrm{zer}_a\, \mathcal{M}$ and closes all $a$-line paths in the graph with edge matrix $\mathcal{M}$ [see (3.1)].

Partial inversion of $M$, with respect to all indices $V = \{1, \ldots, d\}$, gives the inverse of $M$ and partial closure of $\mathcal{M}$, with respect to $V$, gives the edge matrix $\mathcal{M}^-$ of the transitive closure of the graph with edge matrix $\mathcal{M}$

$$
\mathrm{inv}_V\, M = M^{-1}, \qquad \mathrm{zer}_V\, \mathcal{M} = \mathcal{M}^-
$$

[see also equation (2.7)].

Both operators are commutative; that is, for three disjoint subsets $a$, $b$ and $c$ of $V$

$$\operatorname{inv}_a \operatorname{inv}_b M = \operatorname{inv}_b \operatorname{inv}_a M, \qquad \operatorname{zer}_a \operatorname{zer}_b \mathcal{M} = \operatorname{zer}_b \operatorname{zer}_a \mathcal{M},$$

but partial inversion can be undone while partial closure cannot

$$\operatorname{inv}_{ab} \operatorname{inv}_{bc} M = \operatorname{inv}_{ac} M, \qquad \operatorname{zer}_{ab} \operatorname{zer}_{bc} \mathcal{M} = \operatorname{zer}_{abc} \mathcal{M}.$$

For $V = \{a, b\}$, $M$ partially inverted on $a$ coincides with $M^{-1}$ partially inverted on $b$,

$$(A.3) \qquad\qquad\qquad \operatorname{inv}_a M = \operatorname{inv}_b M^{-1}.$$

**A.2. Zero partial regression coefficients and independence.** A Gaussian random vector $Y$ has a nondegenerate distribution if its covariance matrix $\Sigma$ is positive definite. From $\Sigma$ partitioned as $(a, b, c, d)$, the conditional covariance matrix $\Sigma_{ab|c}$, of $Y_a$ and $Y_b$ given $Y_c$, is obtained by partially inverting $\Sigma$, with respect to $c$,

$$\Sigma_{ab|c} = \Sigma_{ab} - \Sigma_{ac} \Sigma_{cc}^{-1} \Sigma_{cb}.$$

The least-squares linear predictor of $Y_a$ given $Y_b$ is $\Pi_{a|b} Y_b$ with $\Pi_{a|b} = \Sigma_{ab} \Sigma_{bb}^{-1}$. For prediction of $Y_a$ with $Y_b$ and $Y_c$, the matrix of the least-squares regression coefficients $\Pi_{a|bc}$ is partitioned as

$$(A.4) \qquad \Pi_{a|bc} = \begin{pmatrix} \Pi_{a|b.c} & \Pi_{a|c.b} \end{pmatrix} = \begin{pmatrix} \Sigma_{ab|c} \Sigma_{bb|c}^{-1} & \Sigma_{ac|b} \Sigma_{cc|b}^{-1} \end{pmatrix},$$

where, for example, $\Pi_{a|b.c} Y_b$ predicts $Y_a$ with $Y_b$ when both $Y_a$ and $Y_b$ are adjusted for linear dependence on $Y_c$. Equation (A.4) is generalized by

$$(A.5) \qquad \Pi_{a|bc.d} = \begin{pmatrix} \Pi_{a|b.cd} & \Pi_{a|c.bd} \end{pmatrix} = \begin{pmatrix} \Sigma_{ab|cd} \Sigma_{bb|cd}^{-1} & \Sigma_{ac|bd} \Sigma_{cc|bd}^{-1} \end{pmatrix}.$$

By property (A.3) of partial inversion

$$\begin{pmatrix} \Pi_{a|b.cd} & \Pi_{a|c.bd} \end{pmatrix} = -(\Sigma^{aa})^{-1} \begin{pmatrix} \Sigma^{ab} & \Sigma^{ac} \end{pmatrix},$$

where $\Sigma^{aa}$, $\Sigma^{ab}$ and $\Sigma^{ac}$ are submatrices of the concentration matrix $\Sigma^{-1}$. From $\Sigma^{-1}$, the concentration matrix $\Sigma^{bc.a}$, of $Y_b$ and $Y_c$ after marginalizing over $Y_a$ is obtained by partially inverting $\Sigma^{-1}$, with respect to $a$,

$$\Sigma^{bc.a} = \Sigma^{bc} - \Sigma^{ba} (\Sigma^{aa})^{-1} \Sigma^{ac}.$$

A recursive relation for matrices of least-squares regression coefficients generalizes a result due to Cochran [1],

$$(A.6) \qquad\qquad \Pi_{a|bc.d} = \Pi_{a|b.c} - \Pi_{a|d.bc} \Pi_{d|b.c},$$

and is obtained with partial inversion of $\Sigma$ first, with respect to $b$, $c$ and $d$.

LEMMA A.1.   *For nondegenerate Gaussian distributions, linear and probabilistic independencies combine equivalently as follows:*

(i) *symmetry:* $\Pi_{a|b.c} = 0$ *implies* $\Pi_{b|a.c} = 0 \iff a \perp\!\!\!\perp b|c$ *implies* $b \perp\!\!\!\perp a|c$;

(ii) *decomposition:* $\Pi_{a|bc.d} = 0$ *implies* $\Pi_{a|b.d} = 0 \iff a \perp\!\!\!\perp bc|d$ *implies* $a \perp\!\!\!\perp b|d$;

(iii) *weak union:* $\Pi_{a|bc.d} = 0$ *implies* $\Pi_{a|b.cd} = 0 \iff a \perp\!\!\!\perp bc|d$ *implies* $a \perp\!\!\!\perp b|cd$;

(iv) *contraction:* $\Pi_{a|b.c} = 0$ *and* $\Pi_{a|d.bc} = 0$ *imply* $\Pi_{a|bd.c} = 0 \iff a \perp\!\!\!\perp b|c$ *and* $a \perp\!\!\!\perp d|bc$ *imply* $a \perp\!\!\!\perp bd|c$;

(v) *intersection:* $\Pi_{a|b.cd} = 0$ *and* $\Pi_{a|c.bd} = 0$ *imply* $\Pi_{a|bc.d} = 0 \iff a \perp\!\!\!\perp b|cd$ *and* $a \perp\!\!\!\perp c|bd$ *imply* $a \perp\!\!\!\perp bc|d$;

(vi) *composition:* $\Pi_{a|c.d} = 0$ *and* $\Pi_{b|c.d} = 0$ *imply* $\Pi_{ab|c.d} = 0 \iff a \perp\!\!\!\perp c|d$ *and* $b \perp\!\!\!\perp c|d$ *imply* $ab \perp\!\!\!\perp c|d$.

PROOF.   Definition (A.4) implies that $\Pi_{a|b.c}$ vanishes if and only if $\Sigma_{ab|c} = 0$, so that (i) results by the symmetry of the conditional covariance matrix. Property (ii) follows by noting that $\Pi_{a|bc.d} = 0$ is equivalent to the vanishing of both $\Sigma_{ab|d}$ and $\Sigma_{ac|d}$ and thus also of $\Pi_{a|b.d} = \Sigma_{ab|d}\Sigma_{bb|d}^{-1}$. Properties (iii) and (v) are direct consequences of (A.5), while (iv) follows with equation (A.6). Finally, (vi) is another consequence of the definition (A.4) and the equality $(\Pi_{a|b})_{\alpha,\beta} = \Pi_{\alpha|\beta.C}$.

The proof is completed by the equivalence of linear and probabilistic independence statements in Gaussian distributions.   □

**Acknowledgments.**   We thank D. R. Cox, K. Sadeghi and a referee for helpful comments, SAMSI, Raleigh, for having arranged a topic-related workshop in 2006.

DIPARTIMENTO DI STATISTICA "G. PARENTI"
UNIVERSITY OF FLORENCE
VIALE MORGAGNI 59
50134 FIRENZE
ITALY
E-MAIL: giovanni.marchetti@ds.unifi.it

MATHEMATICAL STATISTICS
CHALMERS/GÖTEBORGS UNIVERSITET
CHALMERS TVÄRGATA 3
41296 GÖTEBORG
SWEDEN
E-MAIL: wermuth@math.chalmers.se